\documentclass{birkjour}
\usepackage{graphicx}
\usepackage{amsmath} 
\usepackage{psfrag}
\usepackage{enumerate}
\newtheorem{lemma}{Lemma}

 \newtheorem{thm}{Theorem}[section]
 \newtheorem{cor}[thm]{Corollary}

 \theoremstyle{definition}
 
 \newtheorem{algorithm}[lemma]{Algorithm}
 \theoremstyle{remark}

 \numberwithin{equation}{section}

\newcommand{\N}{\mathbb{N}}

\newcommand{\cal}{\mathcal}
\newcommand{\om}{\Omega}
\newcommand{\omc}{{\Omega_c}}

\begin{document}

\title[Parareal intermediate targets methods for optimal control problem]
{Parareal in time intermediate targets methods for optimal control problem}
\thanks{Part of this work was funded by by the ANR-06-CIS6-007-01 project PITAC}
\author[Y. Maday]{Yvon Maday}
\address{UPMC Univ Paris 06, UMR 7598, Laboratoire Jacques-Louis Lions, F-75005, Paris, France,
  and Division of Applied Mathematics, Brown University, Providence, RI, USA.}
\email{maday@ann.jussieu.fr}
\author[M.-K. Riahi]{Mohamed-Kamel Riahi}
\address{UPMC Univ Paris 06, UMR 7598, Laboratoire Jacques-Louis Lions, F-75005, Paris, France.}
\email{riahi@ann.jussieu.fr}
\author[J. Salomon]{Julien Salomon}
\address{CEREMADE, Universit{\'e} Paris-Dauphine, Pl. du Mal. de
  Lattre de Tassigny, F-75016, Paris, France.}
\email{salomon@ceremade.dauphine.fr}



\date{January 15, 2012}

\begin{abstract}
In this paper, we present a method that enables solving in parallel the
  Euler-Lagrange system associated with the optimal control of a parabolic equation.  
Our approach is based on an iterative update of a sequence of intermediate targets that gives rise to independent sub-problems that can be solved in parallel. This method can be coupled with the parareal in time algorithm. Numerical experiments show the efficiency of our method.
\end{abstract}

\maketitle
\section{Introduction}
In the last decade, parallelism across the time \cite{Bur}, based on the decomposition of the time domain has been exploited to
accelerate the simulation of systems governed by time dependent
partial differential equations~\cite{L}.  Among others,
the parareal algorithm~\cite{LMT} or multi-shooting
schemes~\cite{BZ} have shown excellent results. In the framework of
optimal control, this approach has been used to control parabolic
 systems~\cite{MT, MSS}.

  In this paper, we introduce a new approach
 to tackle such problems. The strategy we follow is based on the concept of
 target trajectory that has been introduced in the case of hyperbolic
 systems in~\cite{MST}. Because of the irreversibility of parabolic
 equations, a new definition of this trajectory is considered.
 It enables us to define at each end point of the
 time sub-domains relevant initial conditions and intermediate
 targets, so that the initial problem is split up into independent
 optimization problems.

The paper is organized as follows: the optimal control problem is
introduced in Section~\ref{sec:2} and the parallelization setting is
described in Section~\ref{sec:3}. The properties of the cost
functionals involved in the control problem are studied in
Section~\ref{sec:4}. The general structure of our algorithm is given
in Section~\ref{sec:5} and its convergences is proven in
Section~\ref{sec:6}. In Section~\ref{sec:7}, we propose a fully
parallelized version of our algorithm. Some numerical tests showing
the efficiency of our approach are presented in Section~\ref{sec:8}.

In the sequel, we consider the optimal control problem associated with
the heat equation on a compact  set $\Omega$ and a time interval
$[0,T]$, with $T>0$.
We denote by $\|.\|_\om$ the space norm associated with $L^2(\Omega)$, and
by $\|. \|_\omc$ the
$L^2$-norm corresponding to a sub-domain $\Omega_c\subset \Omega$. Also, we use the notations $\|.\|_v$
(resp. $\|.\|_{v_n}$) and
$\langle . ,. \rangle_v$ (resp.  $\langle . ,. \rangle_{v_n}$) to represent the norm and the scalar product
of the Hilbert space $L^2(0,T;\Omega_c)$ (resp. $L^2(I';\Omega_c)$),
with $I'$ a sub-interval of $[0,T]$ ).  
Given a function $y$ defined on the time interval
$[0,T]$, we denote by $y_{|I'}$ the
restriction of $y$ to $I'$.


\section{Optimal control problem}\label{sec:2}
Given $\alpha>0$, consider the optimal control problem
defined by:
$$\min_{v\in L^2([0,T];L^2(\Omega_c))}J(v),$$
with
$$ J(v)=\frac12 \|y(T)-y_{target}\|_\om^2 +
\dfrac{\alpha}{2}\int_{0}^{T}\|v(t)\|_\omc^2 dt,$$ where $y_{target}$ is a given state in $L^2(\Omega)$.
The state $y$ evolves from $y_0$ on $[0,T]$ according to
$$\partial_t{y} - \nu \Delta y = \cal Bv.$$
In this equation, $\Delta$ denotes the Laplace operator, $v$ is the
control term, applied on $\Omega_c$ and $\cal B$ is
the natural injection from $\Omega_c$ into $\Omega$. We assume Dirichlet
conditions for $y$ on the boundary of $\Omega$.\\
The corresponding optimality system reads as

\begin{equation}\label{os1}
\left\{\begin{array}{ccl}
\partial_t{y} - \nu \Delta y&=& \cal Bv  \qquad \hbox{on} \ [0,T]\times\Omega\\
y(0)&=&y_0,
\end{array}
\right.
\end{equation}

\begin{equation}\label{os2}
\left\{\begin{array}{ccl}
\partial_t{p} + \nu \Delta p &=& 0 \qquad \hbox{on} \ [0,T]\times\Omega\\
p(T) &=&  y(T)-y_{target},
\end{array}
\right.
\end{equation}

\begin{equation}\label{os3}
\alpha v + \cal B^*p =0,
\end{equation}
where $\cal B^*$ is the adjoint operator of $\cal B$.\\

Note that for any $\alpha>0$, the functional $J$ is continuous, $\alpha$-convex
in $L^2(\Omega_c)$ and consequently the system (\ref{os1}--\ref{os3}) has a unique
solution by $v^\star$. We denote by $y^\star$, $p^\star$ the
associated state and adjoint state.  


\section{Time parallelization setting}\label{sec:3}
In this section, we describe the relevant setting for a time parallelized
resolution of the optimality system.\\
Consider $N\geq 1$ and a subdivision of $[0,T]$ of the form:
$$ [0,T]=\cup_{n=0}^{N-1}I_n,$$
with $I_n=[t_n,t_{n+1}]$, $t_0=0<t_1<...<t_{N-1}<t_N=T$. For the sake on simplicity, we assume here that the subdivision is uniform, i.e. 
for $n=0,\dots, N-1$ we assume that
$t_{n+1}-t_n = T/N$
; we denote
$\Delta T=T/N$. Given a control $v$ and its corresponding
state $y$ and adjoint state $p$, we define the {\it target
trajectory} by:

\begin{equation}\label{target_trajectory}
\chi = y - p \qquad \hbox{on} \ [0,T]\times\Omega . 
\end{equation}
The trajectory $\chi$ is not governed by a partial differential
equation, but reaches $\chi(T)  = y_{target}$ at time $T$ from (\ref{os2}$_b$), hence its denomination.\\

For $n=0,\dots,N-1$, consider the sub-problems
\begin{equation}\label{subpb1}
\min_{v_n\in L^2(I_n;L^2(\Omega_c))}J_n(v_n),
\end{equation}
with
\begin{equation}\label{subpb2}
 J_n(v_n)=\frac12 \|y_n(t_{n+1})-\chi(t_{n+1})\|_\om^2 +
\dfrac{\alpha}{2}\int_{I_n}\|v_n(t)\|_\omc^2 dt,
\end{equation}
where the function $y_n$ is defined by
\begin{equation}\label{pos1}
\left\{\begin{array}{ccl}
\partial_t{y_n} - \nu \Delta y_n&=& \cal Bv_n  \qquad \hbox{on} \ I_n\times\Omega\\
y_n(t_n)&=&y(t_n).
\end{array}
\right.
\end{equation}
Recall that this optimal control problem is parameterized by $v$ (and
$y$ and $p$) through the local target  $\chi(t_{n+1})$, we note that this sub-problem has the same structure as the original one, and
is also strictly convex.
The optimality system associated with this
optimization problem is given by \eqref{pos1} and the equations

\begin{equation}\label{pos2}
\left\{\begin{array}{ccl}
\partial_t{p_n} + \nu \Delta p_n &=& 0 \qquad \hbox{on} \ I_n\times\Omega\\
p_n(t_{n+1}) &=&  y(t_{n+1})-\chi(t_{n+1}),
\end{array}
\right.
\end{equation}

\begin{equation}\label{pos3}
\alpha v_n + \cal B^*p_n =0,
\end{equation}
we denote by $v_n^\star$ its solution.


\section{Some properties of $J$ and $J_n$}\label{sec:4}
The introduction of the target trajectory in the last section  is motivated by the following
result.
\begin{lemma}
Denote by $\chi^\star$ the target trajectory defined by
\eqref{target_trajectory} with $y=y^\star$ and $p=p^\star$
and  by $y_n^\star,p_n^\star,v_n^\star$ the solutions of
(\ref{pos1}--\ref{pos3}) with $y=y^\star$ and $\chi=\chi^\star$. One has:
$$ v_n^\star = v^\star_{|I_n}.$$
\end{lemma}
\proof
Thanks to the uniqueness of the solution of the sub-problem, it is
enough to show that $v^\star_{|I_n}$ satisfies the optimality
system~(\ref{pos1}--\ref{pos3}). \\
First, note that $y_{|I_n}^\star$ obviously satisfies~\eqref{pos1}
with $v_n=v^\star_{|I_n}$. It directly follows from the definition of $\chi^\star$
(see~\eqref{target_trajectory}), that:
$$ p^\star(t_{n+1}) =  y^{\star}(t_{n+1})-\chi^\star(t_{n+1}),$$
so that $p^\star_{|I_n}$ satisfies~\eqref{pos2}.
Finally, Equation~\eqref{pos3} is a consequence of~\eqref{os3}. The result follows.
$\hfill \square$

Let $HJ$ denote the hessian operator associated with $J$; there exists a strong connection between the hessian
operators $HJ$ and $HJ_n$ of $J$ and $J_n$, as indicated in the next lemma. 
\begin{lemma}\label{restric}
The hessian operator $HJ_n$ coincides with the restriction of $HJ$ to controls
whose time supports are included in $[t_{N-1},T]$.   
\end{lemma}
\proof
First note that $J$ is quadratic so that $HJ$ is a constant operator.
Given an increase $\delta v\in L^2([0,T];L^2(\Omega_c))$, we have:
$$ \langle HJ(\delta v) , \delta v \rangle_v=\|\delta y (T)\|_\om^2+\alpha\int_0^T\|\delta v(t)\|_\omc^2dt,$$
where $\delta y$ is the solution of
\begin{equation}\label{osdy1}
\left\{\begin{array}{ccl}
\partial_t{\delta y} - \nu \Delta \delta y&=& \cal B \delta v  \qquad
\hbox{on} \ [0,T] \times\Omega\\
\delta y(0)&=&0.
\end{array}
\right.
\end{equation}
Given $1\leq n \leq N$, consider now an increase  $\delta v_n\in L^2(I_n;L^2(\Omega_c))$. One
finds in the same way that:
$$ \langle HJ_n(\delta v_n),   \delta v_n \rangle_{v_n}=\|\delta y_n
(t_{n+1})\|_\om^2+\alpha\int_{t_n}^{t_{n+1}}\|\delta v_n(t)\|^2_\omc dt,$$
where $\delta y_n$ is the solution of
\begin{equation}\label{incr_eq}
\left\{\begin{array}{ccl}
\partial_t{\delta y_n} - \nu \Delta \delta y_n&=& \cal B \delta v_n  \qquad
\hbox{on} \ [t_n,t_{n+1}] \times\Omega\\
\delta y_n(t_n)&=&0.
\end{array}
\right.
\end{equation}
Suppose now that $\delta v=0$ on $[0,t_{N-1}]$, it is a simple 
matter to check that $\delta y\equiv 0$ over $[0, t_{N-1}]$. The restriction of
$\delta y$ on the interval $[t_{N-1},T]$ thus satisfies $\delta
y(t_{N-1})=0$ and is consequently (up to a time translation) the
solution of~\eqref{incr_eq}. $\hfill \square$

We end this section with an estimate on these hessian operators.
\begin{lemma}\label{lemm:poinc}
Given $\delta v\in L^2([0,T];L^2(\Omega_c))$, one has:
\begin{equation}\label{estimhess}
\alpha \int_0^T\|\delta v(t)\|^2_\omc dt \leq\langle   HJ(\delta v) , \delta
v \rangle_v\leq \beta \int_0^T\|\delta v(t)\|^2_\omc dt,
\end{equation}
where $\beta=\alpha + C/\sqrt 2,$ with $C$ the Poincar\'e's constant
associated with $L^2(\Omega)$.
\end{lemma}
The proof of this result is standard and given in Appendix for the sake of completeness. Because of
Lemma~\ref{restric}, the hessian operator $HJ_n$ also
satisfies~\eqref{estimhess}.


\section{Algorithm}\label{sec:5}
We are now in a position to propose a time parallelized procedure
to solve (\ref{os1}--\ref{os3}). In what follows we describe the principal steps of a parallel algorithm named ``{\sc sitpoc}'' ( {\sc s}erial {\sc i}ntermediate {\sc t}argets for {\sc p}arallel {\sc o}ptimal {\sc c}ontrol).
\begin{algorithm}[{\sc sitpoc}]\label{galalg}
Consider an initial control $v^0$ and suppose that, at step $k$ one
knows $v^k$. The computation of $v^{k+1}$ is achieved as follows:
\begin{enumerate}[I.]
\item\label{step1} Compute $y^k$, $p^k$ and the associated target trajectory $\chi^k$ according to
  \eqref{os1}, \eqref{os2} and \eqref{target_trajectory} respectively.
\item\label{step2} Solve approximately the $N$ sub-problems \eqref{subpb1} in
  parallel. For $n=0,\dots,N-1$, denote by $\tilde{v}_n^{k+1}$ the
  corresponding solutions and by $\tilde v^{k+1}$ the concatenation of $(\tilde{v}_n^{k+1})_{n=0,\dots,N-1}$.
\item\label{step4} Define $v^{k+1}$ by $v^{k+1}=(1-\theta^k)v^k+\theta^k\tilde
  v^{k+1}$, where $\theta^k$ is defined to minimize $J((1-\theta^k)v^k+\theta^k\tilde
  v^{k+1})$.
\end{enumerate}
\end{algorithm}
Note that we do not explain in detail here the optimization step (Step \ref{step2}) and
rather present a general  structure of our algorithm.
 Because of the strictly convex setting, some steps of, e.g., a gradient method
 or a small number of conjugate gradient method step can be used. 


\section{Convergence}\label{sec:6}
The convergence of Algorithm~\ref{galalg} can be guaranteed under some
assumptions. In what follows, we denote by $\nabla J$ the gradient of $J$.
\begin{thm}
Suppose  that the sequence $(v^k)_{k\in\N}$ defined in
Algorithm~\ref{galalg} satisfies, for all $k \geq 0$:
\begin{equation}\label{eq0}
J(v^k)\neq J(v^\infty),
\end{equation}  

\begin{equation}\label{eq1}
\langle \nabla J(v^k),v^{k+1} - v^k \rangle_v \leq 0,
\end{equation}  
and 
\begin{equation}\label{eq2}
\|\nabla J(v^k)\|_v\leq \eta\|v^{k+1} - v^{k}\|_v,
\end{equation}
for a given $\eta>0$. Then $(v^k)_{k\in\N}$ converges linearly with
a rate $(1-\frac{2\alpha^2}{\eta^2})\in[0,1)$
to the
solution of~(\ref{os1}--\ref{os3}) 
\end{thm}
Note that in the case~\eqref{eq0} is not satisfied, there exists
$k_0\in\N$ such that $v^{k_0}=v^\infty$ and the optimum is reached in a
finite number of steps.
\proof Define the shifted functional
$$\widetilde{J}(v)=J(v)-J(v^\star),$$
and note that because of the definition of $v^\star$, one has
\begin{equation}\label{loja1}
\widetilde{J}(v)=\frac{1}{2}\langle HJ(v-v^\star),v-v^\star \rangle_v \leq \dfrac\beta 2 \|v-v^\star\|_v^2.
\end{equation} 
 Since $J$ is quadratic, for any $v\in L^2(\Omega_c)$
 $$\nabla J(v)=HJ(v-v^\star), $$
and consequently
 $$\langle \nabla J(v),v-v^\star\rangle_v =\langle
HJ(v-v^\star),v-v^\star\rangle_v\geq \alpha
\|v-v^\star\|_v^2,$$
so that
\begin{equation}\label{loja2}
\|v-v^\star\|_v\leq \frac{1}{\alpha} \|\nabla J(v)\|_v.
\end{equation}
Combining \eqref{loja1} and \eqref{loja2}, one gets
\begin{equation}
\forall v\in L^2(\Omega_c), \qquad \sqrt{\widetilde{J}(v)}\leq \gamma \|\nabla
J(v)\|_v,\label{et2}
\end{equation}
with $\gamma=\frac1{2\sqrt\alpha}$. \\
On the other hand, the variations in the functional between two
iterations of our algorithm reads as
\begin{eqnarray}
J(v^k)-J(v^{k+1})&=&\langle \nabla J(v^k),v^k - v^{k+1} \rangle_v + \frac{1}{2}
\langle  HJ(v^k - v^{k+1})  ,  v^k - v^{k+1}  \rangle_v \nonumber \\
&\geq& \langle \nabla J(v^k),v^k - v^{k+1} \rangle_v + \frac{\alpha}{2} \|v^k -v^{k+1} \|^2_v.\nonumber
\end{eqnarray}
Combining this last inequality with \eqref{eq1}, one finds that :
\begin{equation}\label{et1}
J(v^k)-J(v^{k+1})\geq \dfrac\alpha 2 \|v^k - v^{k+1} \|^2_v .
\end{equation}
Since $\widetilde J(v^k)-\widetilde J(v^{k+1}) = J(v^k)-J(v^{k+1}) \ge 0$, we have:
\begin{eqnarray}
\sqrt{\widetilde{J}(v^k)}-\sqrt{\widetilde{J}(v^{k+1})}&\geq&\dfrac{1}{2\sqrt{\widetilde{J}(v^k)}}\left(J(v^k)-J(v^{k+1})
\right)  \nonumber \\
&\geq& \dfrac{\alpha}{4\sqrt{\widetilde{J}(v^k)}} \|v^k - v^{k+1}
\|^2_v \label{cv1} \\
&\geq& \dfrac{\alpha}{4\gamma \|\nabla
J(v)\|_v} \|v^k - v^{k+1} \|^2_v \label{cv2} \\
&\geq& \dfrac{\alpha}{4\gamma \eta\| v^k - v^{k+1}\|_v} \|v^k -
v^{k+1} \|^2_v \label{cv3} \\
&\geq& c \|v^k - v^{k+1} \|_v \label{cv4} ,
\end{eqnarray}
where $c=\frac{\alpha}{2\gamma \eta}=\frac{\alpha^{\frac32}}\eta$. Indeed \eqref{cv1} follows
from \eqref{et1}, \eqref{cv2} from \eqref{et2} and \eqref{cv3} from \eqref{eq2}.
It follows from the monotonic convergence of $ \sqrt{\widetilde J}(v^{k})$ that the sequence $v^k$ is Cauchy, thus its 
 convergence.
\\
Let us now study the convergence rate. Define
$r^k=\sum_{\ell=k}^{+\infty}\|v^{\ell+1}-v^\ell\|_v$. Summing~\eqref{cv4} between $k$ and $+\infty$, we obtain:
$$ \sqrt{\widetilde{J}(v^{k})}\geq cr^k.$$
Using again~\eqref{et2} and~\eqref{eq2}, one finds that:
\begin{equation}\label{lasteq}
\eta\gamma(r^k-r^{k+1})\geq c r^k.
\end{equation}
Note that this inequality implies that $1-\frac c{\eta\gamma}\geq
0$. Define $C:=\frac{2\alpha^2}{\eta^2}=\frac c{\eta\gamma}$, we
 have $0< C\leq 1$. Because of~\eqref{lasteq}:
$$(1-C)^{-k}r^k\geq(1-C)^{-(k+1)}r^{k+1},$$
and the result follows.
$\hfill
\square$ 
\\\\

We now give an example where hypothesis~(\ref{eq1}--\ref{eq2}) are satisfied.
\begin{cor}~\label{locoptgrad}
Assume that Step~\ref{step2} of Algorithm~\ref{galalg} is achieved using only one step
of a locally optimal step gradient method and that at step $k$, the
algorithm is initialized with $v^k_n:=v^k_{|I_n}$, then~(\ref{eq1}--\ref{eq2}) are 
satisfied hence the algorithm converges to the solution of~(\ref{os1}--\ref{os3}).
\end{cor}
\proof
Because of the assumptions, the optimization step (Step.~\ref{step2}) reads:
$$\tilde v_n^{k+1}=v_n^k-\rho_n^k\nabla J_n(v_n^k).$$
Since the  functionals $J_n$ are quadratic, one has:
$$
\rho_n^k=\frac{\|\nabla J_n(v_n^k)\|^2_{v_n}}{\langle HJ_n(\nabla
  J_n(v_n^k)) , \nabla J_n(v_n^k) \rangle_{v_n}},
$$
A first consequence of these equalities is that:
\begin{equation}\label{ineq1_part}
\langle \nabla J_n(v_n^k),\tilde
v_n^{k+1}-v_n^k\rangle_{v_n}=-\rho_n^k\|\nabla
J_n(v_n^k)\|_{v_n}^2\leq 0.
\end{equation}
Moreover Lemmas~\ref{restric} and~\ref{lemm:poinc} imply:
\begin{equation}\label{rho_estim}
\frac 1\beta \leq \rho_n^k \leq \frac 1\alpha.
\end{equation}
One can also obtain  similar estimates of $\theta^k$. In this view,
note first that since the only iteration which is considered uses as   
directions of descent $\nabla J_n(v^k_n)=\nabla J(v^k)_{|I_n}$.
Then:
\begin{eqnarray*}
\theta^k&=&-\frac{\langle \nabla J(v^k),\tilde v^{k+1}-v^k\rangle_v}{\langle  HJ (\tilde v^{k+1}-v^k)  , \tilde
  v^{k+1}-v^k\rangle_v},\\
&=&-\frac{1}{\langle HJ (\tilde v^{k+1}-v^k) ,  \tilde
  v^{k+1}-v^k\rangle_v}\sum_{n=1}^{N}\frac1{\rho^k_n}\|\tilde v_n^{k+1}-v_n^k\|^2_{v_n}.
\end{eqnarray*}
Using~\eqref{rho_estim}, one deduces:
\begin{equation}
\frac\alpha\beta\leq \theta^k \leq \frac\beta\alpha.
\end{equation}
This preliminary results will now be used to prove the theorem. 
The proof of~\eqref{eq1}, follows from~\eqref{ineq1_part}:
\begin{eqnarray*}
\langle \nabla J(v^k),v^{k+1} - v^k \rangle_v&=&\theta^k\langle \nabla
J(v^k),\tilde v^{k+1} - v^k \rangle_v,\\
&=&\theta^k\sum_{n=1}^{N}\langle \nabla
J_n(v_n^k),\tilde v_n^{k+1} - v_n^k \rangle_{v_n}\leq 0.
\end{eqnarray*}
This last estimate is a consequence of~\eqref{ineq1_part}.
It remains to prove~\eqref{eq2}. We have:
\begin{eqnarray*}
\|v^{k+1}-v^k\|_v&=&\theta^k\|\tilde v^{k+1}-v^k\|_v\\
&=&\theta^k\sqrt{\sum_{n=1}^{N}\|\tilde v_n^{k+1}-v_n^k\|^2_{v_n}}\\
&=&\theta^k\sqrt{\sum_{n=1}^{N}\left(\rho_n^k\right)^2\|\nabla
  J_n(v_n^k)\|^2_{v_n}}\\
&\leq& \frac\alpha\beta\sqrt{\sum_{n=1}^{N}\frac1{\beta^2}\|\nabla
  J_n(v_n^k)\|^2_{v_n}}\\
&=&\frac{\alpha}{\beta^2}\|\nabla
  J(v^k)\|_{v},
\end{eqnarray*}
and the result follows.

$\hfill
\square$


\section{Parareal acceleration}\label{sec:7}
The method we have presented with algorithm~\ref{galalg} requires in Step~\ref{step1} two sequential resolutions of
the evolution Equation~\eqref{os1} on the whole interval $[0,T]$, which does
not fit with the parallel setting. In this section, we make use of the
parareal algorithm to parallelize the corresponding computations.
\subsection{Setting}
Let us first recall the main features of the parareal
algorithm. We consider the example of Equation~\eqref{os1}. 
In order to solve in parallel an evolution equation, for the parareal
scheme~\cite{L} we introduce intermediate initial conditions at times $(t_n)_{n=0,...,N-1}$ that are
updated iteratively. Suppose that these values $(\lambda_n^k)$ are
known at step $k$. Denote by $\cal G_n(\lambda_n)$ and $\cal
F_n(\lambda_n)$ coarse and fine solutions of~\eqref{pos1} at time
$t_{n+1}$ with $\lambda_n$ as initial value. The update is done according to the following
iteration:
$$\lambda_{n+1}^{k+1}=\cal G_n (\lambda_n^{k+1})+\cal F_n (\lambda_n^k)
-\cal G_n (\lambda_n^k).$$
We use this procedure in Step~\ref{step1} of
Algorithm~\ref{galalg}. The idea we follow consists in merging the two
procedures, i.e. doing one parareal iteration at each iteration of our
algorithm.  
\subsection{Algorithm}
We now give details on the resulting procedure. Since the evolution
equations depend on the control, we replace the notations $\cal
G_n(\lambda_n)$ and $\cal F_n(\lambda_n)$ by $\cal
G_n(\lambda_n,v_n)$ and $\cal F_n(\lambda_n,v_n)$ respectively. As we
 need backward solvers to compute $p$, see~\eqref{os2}, we
also introduce  $\widetilde{\cal G}_n(\mu_{n+1})$ and
$\widetilde{\cal F}_n(\mu_{n+1})$ to denote coarse and fine solutions of~\eqref{pos2} at time
$t_{n}$ with $\mu_{n+1}$ as ``initial'' value (given at time $t_{n+1}$). Note that these bakward solvers $\widetilde{\cal F_{n}}$ (resp: $\widetilde{\cal G_{n}}$) do not depend on the control.\\

We describe in the following the principal steps of an enhanced version of the  {\sc sitpoc} algorithm which we give the name ``{\sc pitpoc}'' as {\sc p}arareal {\sc i}ntermediate {\sc t}argets for {\sc o}ptimal {\sc c}ontrol. 
\begin{algorithm}[{\sc pitpoc}]\label{paralg}
Denote by $v_n^k = v^k_{|I_n}$. Consider a control $(v_{n}^0)_{n=0,...,N-1}$, initial values $(\lambda^0_n)_{n=0,...,N}$ (through forward scheme $\lambda_{n+1}^0=\cal G_{n}(\lambda_{n}^0,v_{n}^0)$),
final values $(\mu^0_{n})_{n=1,\dots,N}$  (through backward scheme $\mu_{n}^{0}=\widetilde{\cal G}_n(\mu_{n+1}^{0})$.\\
 Suppose that, at
step $k$ one  
knows $v^k$, $(\lambda^k_n)_{n=0,...,N}$ and 
$(\mu^k_{n})_{n=1,\dots,N}$. The computation of $v^{k+1}$,
$(\lambda^{k+1}_n)_{n=0,...,N}$ and 
$(\mu^{k+1}_{n})_{n=1,\dots,N}$ is
achieved as follows: 
\begin{enumerate}[I.]
 \item Build the target trajectory $(\chi^{k}_n)_{n=1,...,N}$ according to
a definition similar to~\eqref{target_trajectory}:
$$\chi_{n}^{k} = \lambda_{n}^{k} - \mu_{n}^{k}.$$
\item\label{step2para} Solve approximately the $N$ sub-problems \eqref{subpb1}  in
  parallel. For $n=0,\dots,N-1$, denote by $\tilde{v}_n^{k+1}$ the
  corresponding solutions.
\item Define $\tilde v^{k+1}$ as the concatenation of the sequence $(\tilde{v}_n^{k+1})_{n=0,\dots,N-1}$.
\item\label{step1para} Compute $(\tilde\lambda^{k+1}_n)_{n=0,...,N}$, 
$(\mu^{k+1}_{n})_{n=1,\dots,N}$ by:
\begin{eqnarray*}
\tilde \lambda_{n+1}^{k+1}&=&\cal G_n (\tilde \lambda_n^{k+1},\tilde v_n^{k+1})+\cal F_n (\lambda_n^k,\tilde v_n^{k+1})-\cal G_n (\lambda_n^k,v_n^k),\\
\mu_{n}^{k+1}&=&\widetilde{\cal G}_n(\mu_{n+1}^{k+1})+\widetilde{\cal F}_n (\mu_{n+1}^k)-\widetilde{\cal G}_n (\mu_{n+1}^k), 
\end{eqnarray*}
\item\label{step4para} Define $v^{k+1}$ and
  $(\lambda^{k+1}_n)_{n=0,...,N}$
\begin{eqnarray*}
      v^{k+1}&=&(1-\theta^k)v^k      +\theta^k\tilde v^{k+1},\\
\lambda_n^{k+1}&=&(1-\theta^k)\lambda_n^k+\theta^k\tilde \lambda_n^{k+1}
\end{eqnarray*}
 where $\theta^k$ is defined to minimize 
$$\frac12 \|(1-\theta^k)\lambda_N^k+\theta^k\tilde \lambda_N^{k+1}-y_{target}\|_\om^2 +
\dfrac{\alpha}{2}\int_{0}^{T}\| (1-\theta^k)v^k(t)      +\theta^k\tilde v^{k+1}(t)\|_\omc^2 dt.$$
\item $k=k+1$ and return to I.
\end{enumerate}
\end{algorithm}

\section{Numerical Results}\label{sec:8}
In this section, we test the efficiency of our method and show how robust the approach is. We consider two independent parts describing numerical results of the selected algorithm. 
\subsection{Setting}
We consider a
2D example, where $\Omega=[0,1]\times [0,1]$ and
$\Omega_c=[\frac{1}{3},\frac{2}{3}]\times[\frac{1}{3},\frac{2}{3}]$. The
parameters related to our control problem are $T=6.4$, 
$\alpha=10^{-2}$ and $\nu=10^{-2}$. The time interval is discretized using
a uniform step $\delta t=10^{-2}$, and an Implicit-Euler solver is used
to approximate the solution of Equations~(\ref{os1}--\ref{os2}). For the space
discretization, we use $\mathbb{P}_{1}$ finite elements. Our
implementation makes use of the freeware {\texttt FreeFem}~\cite{bib7}
and the parallelization is achieved thanks to the Message Passing
Interface library. The independent optimization procedures required in
Step~\ref{step2} are simply carried out using one step of an optimal gradient method. 
\subsection{Influence of the number of sub-intervals}
\begin{figure}[h!]\centering
   \begin{minipage}[c]{0.4\linewidth}\centering
       \includegraphics[width=6.3cm,height=6cm]{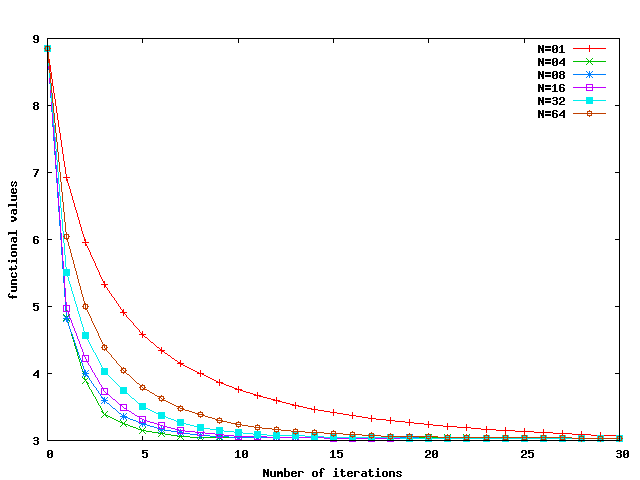}
   \end{minipage}\hfill
   \begin{minipage}[c]{0.4\linewidth} \centering 
      \includegraphics[width=6.3cm,height=6cm]{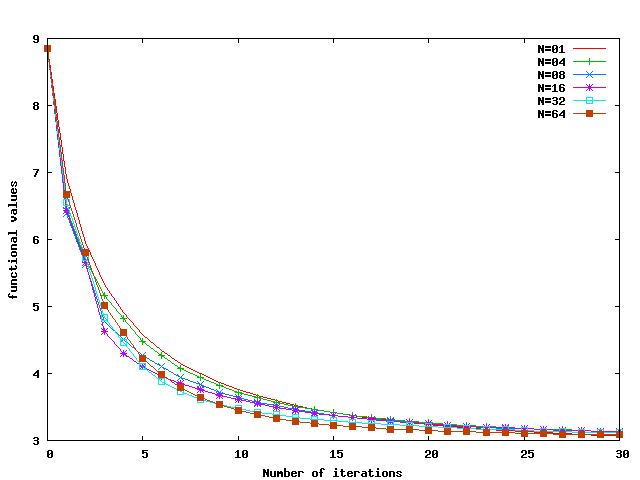}
   \end{minipage}
  \caption{Decaying cost functional values according to the iterations count with respect to {\sc sitpoc} algorithm (left) and {\sc pitpoc} algorithm (right).}
   \label{itersitpit}
\end{figure}

\begin{figure}[hpt!]\centering
   \begin{minipage}[c]{0.4\linewidth}\centering
       \includegraphics[width=6.3cm,height=6cm]{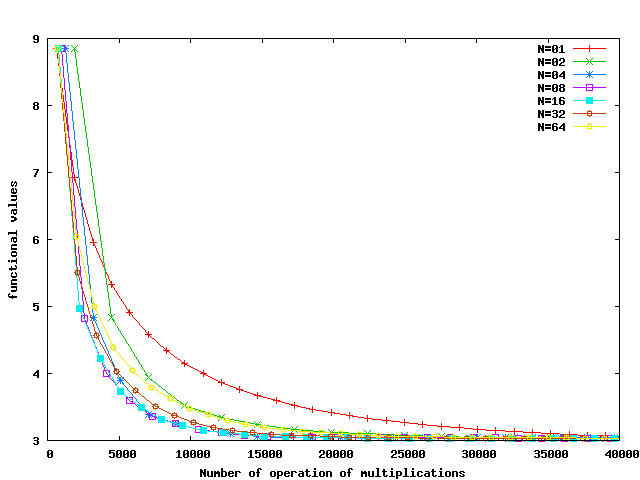}
   \end{minipage}\hfill
   \begin{minipage}[c]{0.4\linewidth} \centering 
      \includegraphics[width=6.3cm,height=6cm]{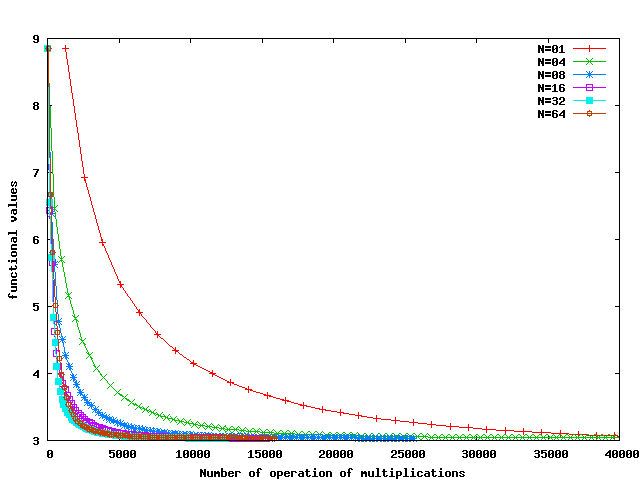}
   \end{minipage}
  \caption{Decaying cost functional values according to the multiplication operations  count with respect to {\sc sitpoc} algorithm (left) and {\sc pitpoc} algorithm (right).}
   \label{multsitpit}
\end{figure}

\begin{figure}[hpt!]\centering
   \begin{minipage}[c]{0.4\linewidth}\centering
       \includegraphics[width=6.3cm,height=6cm]{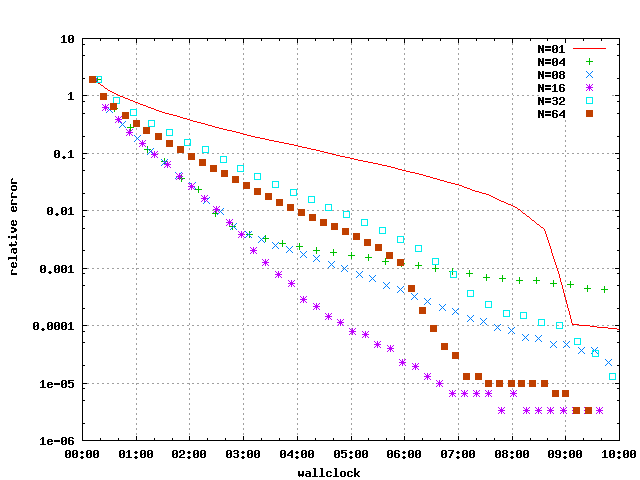}
   \end{minipage}\hfill
   \begin{minipage}[c]{0.4\linewidth} \centering 
      \includegraphics[width=6.3cm,height=6cm]{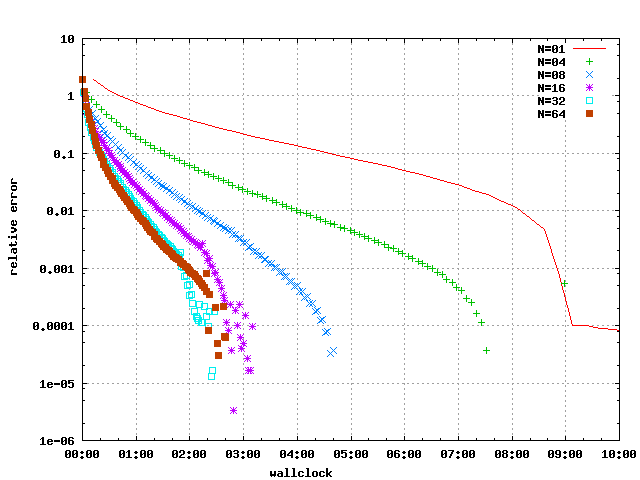}
   \end{minipage}
  \caption{Decaying cost functional values according to elapsed real time with respect to {\sc sitpoc} algorithm (left) and {\sc pitpoc} algorithm (right).}
   \label{timesitpit}
\end{figure}

In this section, Step II of Algorithm~\ref{galalg} and Algorithm~\ref{paralg} are achieved by using one step of an optimal step gradient method. We first test our algorithm by varying the number of sub-intervals. The evolution of the cost functional values are plotted with respect to the number of iteration (Figure~\ref{itersitpit}), the number of matrix multiplication (Figure~\ref{multsitpit}) and the number of wall-clock time of computation (Figure~\ref{timesitpit}).
We first note that Algorithm 4 actually acts as a preconditioner, since it improves the convergence rate of the optimization process. The introduceion of the intermediates targets allows to accelerate the decrease of 
the functional values, as shown in Figure~\ref{multsitpit} (left). Note that this property holds mostly for small numbers of sub-intervals, and disapears when dealing with large subdivisions.
This feature is lost when considering Algorithm~\ref{paralg}, whose convergence does not significantly depend on the number of sub-intervals that is considered, see Figure~\ref{multsitpit} (right). 	

On the contrary, Algorithm~\ref{paralg} achieves a good acceleration when considering the number of mutliplications involved in the computations. The corresponding results are shown in Figure~\ref{multsitpit}, where the parallel operations have been counted only once. We see that Algorithm  is close to the full efficiency, since the number of multiplications required to obtain a given value for the cost functional is roughly proportional to $\frac 1N$. 

We finally consider the wall-clock time required to carry out our algorithms. As the main part of the operations involved in the computation consists in matrix multiplications, the results we present in Figure 3 are close to the ones of Figure 2.

\subsection{Influence of the number of steps in the optimization method}
We now vary the number of steps of the gradient method used in Step II of our algorithm. The results are presented in Figure 4. Subdivisions of $N=4$ and $N=16$ intervals are considered. In both cases, we see that an increase in the number of gradient steps improves the preconditionning feature of our algorithm. However, we also observe that this strategy saturates for large numbers of gradient steps which probably reveals that the sub-problems considered in Step II are practically solved after 5 sub-iterations.  
\medskip

More results can be found in \cite{theseK}.

\begin{figure}[hpt!]\centering
   \begin{minipage}[c]{0.4\linewidth}\centering
       \includegraphics[width=6.3cm,height=6cm]{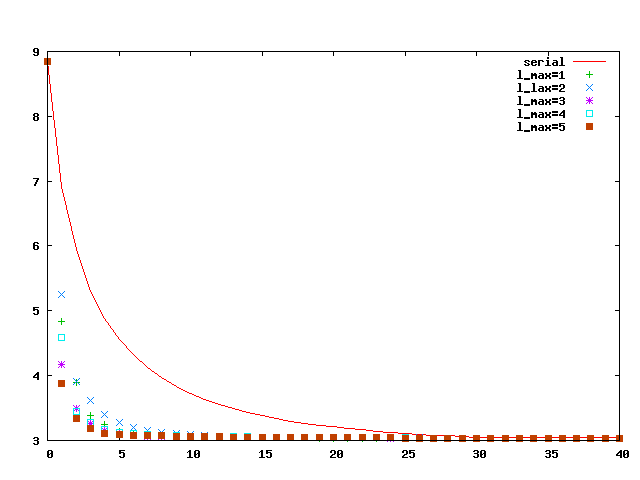}
   \end{minipage}\hfill
   \begin{minipage}[c]{0.4\linewidth} \centering 
      \includegraphics[width=6.3cm,height=6cm]{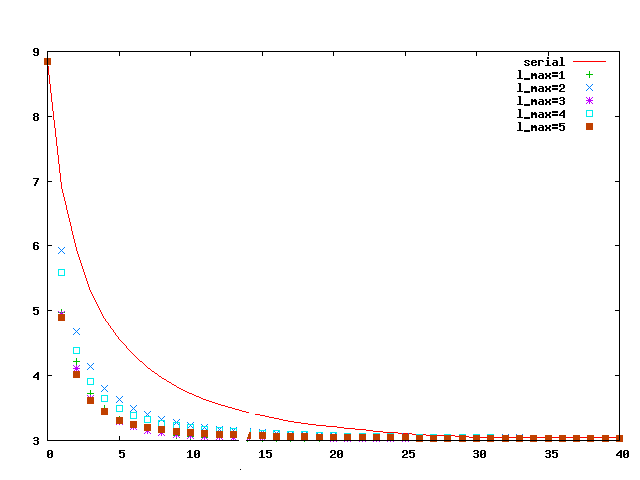}
   \end{minipage}
  \caption{{\sc sitpoc} algorithm with 4 subdivisions (left) and $16$ subdivision (right): variation of the number of (lower/local)inner-iterations $\ell_\text{max}$.}
\label{iterloc0416}
\end{figure}

\section*{Appendix}
For the sake of completeness, we recall here the proof of
Lemma~\ref{lemm:poinc}.\\
Because of~\eqref{osdy1} and thanks to Young's inequality, one has for
all $t\in[0,T]$ and all $\varepsilon>0$:
\begin{eqnarray}
\frac 12 \frac{d}{dt}\|\delta y(t)\|_\om^2 + \nu\|\nabla_x \delta
y(t)\|_\om^2&=&\int_\Omega\delta y(t) \cal B\delta  v(t) dt \nonumber\\
&\leq &\frac 12\left(\varepsilon \|\delta y(t)\|_\om^2+\frac
1\varepsilon \|\cal B\delta v(t)\|_v^2 \right),\label{young}
\end{eqnarray}
where $\nabla_x$ denotes the gradient with respect to the space variable.
As $\delta y$ is supposed to satisfies Dirichlet conditions, one can
apply Poincar\'e's inequality to obtain:
$$ \|\delta y(t)\|_\om\leq C \|\nabla_x \delta
y(t)\|_\om,$$
for a given $C>0$. Combining this last estimate with~\eqref{young},
one gets:
\begin{equation*}
\frac 12 \frac{d}{dt}\|\delta y(t)\|_\om^2\leq \left(\frac \varepsilon 2
- \frac \nu{C^2}\right)\|\delta y(t)\|_\om^2+\frac
1{2\varepsilon} \|\cal B\delta v(t)\|_v^2.
\end{equation*}
Now, setting $\varepsilon=\frac {2\nu} {C^2}$ gives:
\begin{equation*}
 \frac{d}{dt}\|\delta y(t)\|_\om^2\leq \frac 1\varepsilon \|\cal B\delta v(t)\|_v^2.
\end{equation*}
Since $\|\delta y(0)\|_\om^2=0$, the result follows with the fact that $\|\cal B\|_{2}\leq1$.








\end{document}